# EFFECT OF MEAN ON VARIANCE FUNCTION ESTIMATION IN NONPARAMETRIC REGRESSION

By Lie Wang, Lawrence D. Brown, T. Tony Cai[1]
and Michael Levine

*University of Pennsylvania, University of Pennsylvania, University of Pennsylvania and Purdue University*

Variance function estimation in nonparametric regression is considered and the minimax rate of convergence is derived. We are particularly interested in the effect of the unknown mean on the estimation of the variance function. Our results indicate that, contrary to the common practice, it is not desirable to base the estimator of the variance function on the residuals from an optimal estimator of the mean when the mean function is not smooth. Instead it is more desirable to use estimators of the mean with minimal bias. On the other hand, when the mean function is very smooth, our numerical results show that the residual-based method performs better, but not substantial better than the first-order-difference-based estimator. In addition our asymptotic results also correct the optimal rate claimed in Hall and Carroll [*J. Roy. Statist. Soc. Ser. B* **51** (1989) 3–14].

**1. Introduction.** Consider the heteroscedastic nonparametric regression model

(1) $$y_i = f(x_i) + V^{1/2}(x_i)z_i, \qquad i = 1, \ldots, n,$$

where $x_i = i/n$ and $z_i$ are independent with zero mean, unit variance and uniformly bounded fourth moments. Both the mean function $f$ and variance function $V$ are defined on $[0,1]$ and are unknown. The main object of interest is the variance function $V$. The estimation accuracy is measured both globally by the mean integrated squared error

(2) $$R(\hat{V}, V) = E \int_0^1 (\hat{V}(x) - V(x))^2 \, dx$$

Received March 2006; revised April 2007.
[1]Supported in part by NSF Grant DMS-03-06576.
*AMS 2000 subject classifications.* 62G08, 62G20.
*Key words and phrases.* Minimax estimation, nonparametric regression, variance estimation.

and locally by the mean squared error at a point

$$(3) \qquad R(\hat{V}(x_*), V(x_*)) = E(\hat{V}(x_*) - V(x_*))^2.$$

We wish to study the effect of the unknown mean $f$ on the estimation of the variance function $V$. In particular, we are interested in the case where the difficulty in estimation of $V$ is driven by the degree of smoothness of the mean $f$.

The effect of not knowing the mean $f$ on the estimation of $V$ has been studied before in Hall and Carroll (1989). The main conclusion of their paper is that it is possible to characterize explicitly how the smoothness of the unknown mean function influences the rate of convergence of the variance estimator. In association with this they claim an explicit minimax rate of convergence for the variance estimator under pointwise risk. For example, they state that the "classical" rates of convergence $(n^{-4/5})$ for the twice differentiable variance function estimator is achievable if and only if $f$ is in the Lipschitz class of order at least $1/3$. More precisely, Hall and Carroll (1989) stated that, under the pointwise mean squared error loss, the minimax rate of convergence for estimating $V$ is

$$(4) \qquad \max\{n^{-4\alpha/(2\alpha+1)}, n^{-2\beta/(2\beta+1)}\}$$

if $f$ has $\alpha$ derivatives and $V$ has $\beta$ derivatives. We shall show here that this result is in fact incorrect.

In the present paper we revisit the problem in the same setting as in Hall and Carroll (1989). We show that the minimax rate of convergence under both the pointwise squared error and global integrated mean squared error is

$$(5) \qquad \max\{n^{-4\alpha}, n^{-2\beta/(2\beta+1)}\}$$

if $f$ has $\alpha$ derivatives and $V$ has $\beta$ derivatives. The derivation of the minimax lower bound is involved and is based on a moment matching technique and a two-point testing argument. A key step is to study a hypothesis testing problem where the alternative hypothesis is a Gaussian location mixture with a special moment matching property. The minimax upper bound is obtained using kernel smoothing of the squared first order differences.

Our results have two interesting implications. First, if $V$ is known to belong to a regular parametric model, such as the set of positive polynomials of a given order, the cutoff for the smoothness of $f$ on the estimation of $V$ is $1/4$, not $1/2$ as stated in Hall and Carroll (1989). That is, if $f$ has at least $1/4$ derivative then the minimax rate of convergence for estimating $V$ is solely determined by the smoothness of $V$ as if $f$ were known. On the other hand, if $f$ has less than $1/4$ derivative then the minimax rate depends on the relative smoothness of both $f$ and $V$ and will be completely driven by the roughness of $f$.

Second, contrary to the common practice, our results indicate that it is often not desirable to base the estimator $\hat{V}$ of the variance function $V$ on the residuals from an *optimal* estimator $\hat{f}$ of the mean function $f$ when $f$ is not smooth. Instead it is more desirable to use estimators of the mean $f$ with minimal bias. The main reason is that the bias and variance of $\hat{f}$ have quite different effects on the estimation of $V$. The bias of $\hat{f}$ cannot be removed or even reduced in the second stage smoothing of the squared residuals, while the variance of $\hat{f}$ can be incorporated easily. On the other hand, when the mean function is very smooth, our numerical results show that the residual-based method performs better, but not substantial better than the first-order-difference-based estimator.

The paper is organized as follows. Section 2 presents an upper bound for the minimax risk while Section 3 derives a rate-sharp lower bound for the minimax risk under both the global and local losses. The lower and upper bounds together yield the minimax rate of convergence. Section 4 discusses the obtained results and their implications for practical variance estimation in the nonparametric regression. Section 5 considers finite sample performance of the difference-based method for estimating the variance function. The proofs are given in Section 6.

**2. Upper bound.** In this section we shall construct a kernel estimator based on the square of the first order differences. Such and more general difference based kernel estimators of the variance function have been considered, for example, in Müller and Stadtmüller (1987, 1993). For estimating a constant variance, difference based estimators have a long history. See von Neumann (1941, 1942), Rice (1984), Hall, Kay and Titterington (1990) and Munk, Bissantz, Wagner and Freitag (2005).

Define the Lipschitz class $\Lambda^\alpha(M)$ in the usual way,

$$\Lambda^\alpha(M) = \{g : \text{for all } 0 \leq x, y \leq 1, \ k = 0, \ldots, \lfloor \alpha \rfloor - 1,$$
$$|g^{(k)}(x)| \leq M \text{ and } |g^{(\lfloor \alpha \rfloor)}(x) - g^{(\lfloor \alpha \rfloor)}(y)| \leq M|x - y|^{\alpha'}\},$$

where $\lfloor \alpha \rfloor$ is the largest integer less than $\alpha$ and $\alpha' = \alpha - \lfloor \alpha \rfloor$. We shall assume that $f \in \Lambda^\alpha(M_f)$ and $V \in \Lambda^\beta(M_V)$. We say that the function $f$ "has $\alpha$ derivative" if $f \in \Lambda^\alpha(M_f)$ and $V$ "has $\beta$ derivatives" if $V \in \Lambda^\beta(M_V)$.

For $i = 1, 2, \ldots, n - 1$, set $D_i = y_i - y_{i+1}$. Then one can write

(6) $\quad D_i = f(x_i) - f(x_{i+1}) + V^{1/2}(x_i)z_i - V^{1/2}(x_{i+1})z_{i+1} = \delta_i + \sqrt{2}V_i^{1/2}\epsilon_i,$

where $\delta_i = f(x_i) - f(x_{i+1})$, $V_i^{1/2} = \sqrt{\frac{1}{2}(V(x_i) + V(x_{i+1}))}$ and

$$\epsilon_i = (V(x_i) + V(x_{i+1}))^{-1/2}(V^{1/2}(x_i)z_i - V^{1/2}(x_{i+1})z_{i+1})$$

has zero mean and unit variance.

We construct an estimator $\hat{V}$ by applying kernel smoothing to the squared differences $D_i^2$ which have means $\delta_i^2 + 2V_i$. Let $K(x)$ be a kernel function satisfying

$$K(x) \text{ is supported on } [-1,1], \int_{-1}^{1} K(x)\,dx = 1,$$

$$\int_{-1}^{1} K(x)x^i\,dx = 0 \quad \text{for} \quad i=1,2,\ldots,\lfloor\beta\rfloor \quad \text{and}$$

$$\int_{-1}^{1} K^2(x)\,dx = k < \infty.$$

It is well known in kernel regression that special care is needed in order to avoid significant, sometimes dominant, boundary effects. We shall use the boundary kernels with asymmetric support, given in Gasser and Müller (1979, 1984), to control the boundary effects. For any $t \in [0,1]$, there exists a boundary kernel function $K_t(x)$ with support $[-1,t]$ satisfying the same conditions as $K(x)$, that is,

$$\int_{-1}^{t} K_t(x)\,dx = 1,$$

$$\int_{-1}^{t} K_t(x)x^i\,dx = 0 \quad \text{for } i=1,2,\ldots,\lfloor\beta\rfloor,$$

$$\int_{-1}^{t} K_t^2(x)\,dx \leq \widehat{k} < \infty \quad \text{for all } t \in [0,1].$$

We can also make $K_t(x) \to K(x)$ as $t \to 1$ (but this is not necessary here). See Gasser, Müller and Mammitzsch (1985). For any $0 < h < 1/2$, $x \in [0,1]$, and $i = 2,\ldots,n-2$, let

$$K_i^h(x) = \begin{cases} \int_{(x_i+x_{i-1})/2}^{(x_i+x_{i+1})/2} \frac{1}{h} K\left(\frac{x-u}{h}\right) du, & \text{when } x \in (h, 1-h), \\ \int_{(x_i+x_{i-1})/2}^{(x_i+x_{i+1})/2} \frac{1}{h} K_t\left(\frac{x-u}{h}\right) du, & \text{when } x = th \text{ for some } t \in [0,1], \\ \int_{(x_i+x_{i-1})/2}^{(x_i+x_{i+1})/2} \frac{1}{h} K_t\left(-\frac{x-u}{h}\right) du, & \\ & \text{when } x = 1-th \text{ for some } t \in [0,1], \end{cases}$$

and we take the integral from 0 to $(x_1 + x_2)/2$ for $i = 1$, and from $(x_{n-1} + x_{n-2})/2$ to 1 for $i = n-1$. Then we can see that for any $0 \leq x \leq 1$, $\sum_{i=1}^{n-1} K_i^h(x) = 1$. Define the estimator $\widehat{V}$ as

(7) $$\widehat{V}(x) = \tfrac{1}{2} \sum_{i=1}^{n-1} K_i^h(x) D_i^2.$$

Same as in the mean function estimation problem, the optimal bandwidth $h_n$ can be easily seen to be $h_n = O(n^{-1/(1+2\beta)})$ for $V \in \Lambda^\beta(M_V)$. For this optimal choice of the bandwidth, we have the following theorem.

THEOREM 1. *Under the regression model* (1) *where* $x_i = i/n$ *and* $z_i$ *are independent with zero mean, unit variance and uniformly bounded fourth moments, let the estimator* $\widehat{V}$ *be given as in* (7) *with the bandwidth* $h = O(n^{-1/(1+2\beta)})$. *Then there exists some constant* $C_0 > 0$ *depending only on* $\alpha$, $\beta$, $M_f$ *and* $M_V$ *such that for sufficiently large* $n$,

$$\sup_{f \in \Lambda^\alpha(M_f), V \in \Lambda^\beta(M_V)} \sup_{0 \le x_* \le 1} E(\widehat{V}(x_*) - V(x_*))^2$$
(8)
$$\le C_0 \cdot \max\{n^{-4\alpha}, n^{-2\beta/(1+2\beta)}\}$$

*and*

$$\sup_{f \in \Lambda^\alpha(M_f), V \in \Lambda^\beta(M_V)} E \int_0^1 (\widehat{V}(x) - V(x))^2 \, dx$$
(9)
$$\le C_0 \cdot \max\{n^{-4\alpha}, n^{-2\beta/(1+2\beta)}\}.$$

REMARK 1. The uniform rate of convergence given in (8) yields immediately the pointwise rate of convergence that for any fixed point $x_* \in [0,1]$

$$\sup_{f \in \Lambda^\alpha(M_f), V \in \Lambda^\beta(M_V)} E(\widehat{V}(x_*) - V(x_*))^2 \le C_0 \cdot \max\{n^{-4\alpha}, n^{-2\beta/(1+2\beta)}\}.$$

REMARK 2. It is also possible to use the local linear regression estimator instead of the Priestley–Chao kernel estimator. In this case, the boundary adjustment is not necessary as it is well known that the local linear regression adjusts automatically in boundary regions, preserving the asymptotic order of the bias intact. However, the proof is slightly more technically involved when using the local linear regression estimator. For details see, for example, Fan and Gijbels (1996).

REMARK 3. It is important to note here that the results given in Theorem 1 can be easily generalized to the case of random design. In particular, if the observations $X_1, \ldots, X_n$ are i.i.d. with the design density $f(x)$ that is bounded away from zero (i.e., $f(x) \ge \delta > 0$ for all $x \in [0,1]$), then the results of Theorem 1 are still valid conditionally. In other words,

$$\sup_{f \in \Lambda^\alpha(M_f), V \in \Lambda^\beta(M_V)} \sup_{0 \le x_* \le 1} E(\widehat{V}(x_*) - V(x_*))^2 | X_1, \ldots, X_n)$$
$$\le C_0 \cdot \max\{n^{-4\alpha}, n^{-2\beta/(1+2\beta)}\} + o_p(\max\{n^{-4\alpha}, n^{-2\beta/(1+2\beta)}\})$$

and

$$\sup_{f\in\Lambda^\alpha(M_f),V\in\Lambda^\beta(M_V)} E\left(\int_0^1 (\widehat{V}(x) - V(x))^2 \, dx | X_1,\ldots,X_n\right)$$
$$\leq C_0 \cdot \max\{n^{-4\alpha}, n^{-2\beta/(1+2\beta)}\} + o_p(\max\{n^{-4\alpha}, n^{-2\beta/(1+2\beta)}\})$$

where the constant $C_0 > 0$ now also depends on $\delta$.

**3. Lower bound.** In this section we derive a lower bound for the minimax risk of estimating the variance function $V$ under the regression model (1). The lower bound shows that the upper bound given in the previous section is rate-sharp. As in Hall and Carroll (1989) we shall assume in the lower bound argument that the errors are normally distributed, that is, $z_i \overset{\text{i.i.d.}}{\sim} N(0,1)$.

THEOREM 2. *Under the regression model* (1) *with* $z_i \overset{\text{i.i.d.}}{\sim} N(0,1)$,

$$(10) \quad \inf_{\widehat{V}} \sup_{f\in\Lambda^\alpha(M_f),V\in\Lambda^\beta(M_V)} E\|\widehat{V} - V\|_2^2 \geq C_1 \cdot \max\{n^{-4\alpha}, n^{-2\beta/(1+2\beta)}\}$$

*and for any fixed* $x_* \in (0,1)$

$$(11) \quad \inf_{\widehat{V}} \sup_{f\in\Lambda^\alpha(M_f),V\in\Lambda^\beta(M_V)} E(\widehat{V}(x_*) - V(x_*))^2$$
$$\geq C_1 \cdot \max\{n^{-4\alpha}, n^{-2\beta/(1+2\beta)}\},$$

*where* $C_1 > 0$ *is a constant depending only on* $\alpha$, $\beta$, $M_f$ *and* $M_V$.

It follows immediately from Theorems 1 and 2 that the minimax rate of convergence for estimating $V$ under both the global and local losses is

$$\max\{n^{-4\alpha}, n^{-2\beta/(1+2\beta)}\}.$$

The proof of this theorem can be naturally divided into two parts. The first step is to show

$$(12) \quad \inf_{\widehat{V}} \sup_{f\in\Lambda^\alpha(M_f),V\in\Lambda^\beta(M_V)} E(\widehat{V}(x_*) - V(x_*))^2 \geq C_1 n^{-2\beta/(1+2\beta)}.$$

This part is standard and relatively easy. Brown and Levine (2006) contains a detailed proof of this assertion for the case $\beta = 2$. Their argument can be easily generalized to other values of $\beta$. We omit the details.

The proof of the second step,

$$(13) \quad \inf_{\widehat{V}} \sup_{f\in\Lambda^\alpha(M_f),V\in\Lambda^\beta(M_V)} E(\widehat{V}(x_*) - V(x_*))^2 \geq C_1 n^{-4\alpha},$$

is much more involved. The derivation of the lower bound (13) is based on a moment matching technique and a two-point testing argument. One of the main steps is to study a complicated hypothesis testing problem where the alternative hypothesis is a Gaussian location mixture with a special moment matching property.

More specifically, let $X_1, \ldots, X_n \overset{\text{i.i.d.}}{\sim} P$ and consider the following hypothesis testing problem between

$$H_0 : P = P_0 = N(0, 1 + \theta_n^2)$$

and

$$H_1 : P = P_1 = \int N(\theta_n \nu, 1) G(d\nu),$$

where $\theta_n > 0$ is a constant and $G$ is a distribution of the mean $\nu$ with compact support. The distribution $G$ is chosen in such a way that, for some positive integer $q$ depending on $\alpha$, the first $q$ moments of $G$ match exactly with the corresponding moments of the standard normal distribution. The existence of such a distribution is given in the following lemma from Karlin and Studden (1966).

LEMMA 1. *For any fixed positive integer $q$, there exist a $B < \infty$ and a symmetric distribution $G$ on $[-B, B]$ such that $G$ and the standard normal distribution have the same first $q$ moments, that is,*

$$\int_{-B}^{B} x^j G(dx) = \int_{-\infty}^{+\infty} x^j \varphi(x) \, dx, \qquad j = 1, 2, \ldots, q,$$

*where $\varphi$ denotes the density of the standard normal distribution.*

The moment matching property makes the testing between the two hypotheses "difficult." The lower bound (13) then follows from a two-point argument with an appropriately chosen $\theta_n$. Technical details of the proof are given in Section 6.

REMARK 4. For $\alpha$ between $1/4$ and $1/8$, a much simpler proof can be given with a two-point mixture for $P_1$ which matches the mean and variance, but not the higher moments, of $P_0$ and $P_1$. However, this simpler proof fails for smaller $\alpha$. It appears to be necessary in general to match higher moments of $P_0$ and $P_1$.

REMARK 5. Hall and Carroll (1989) gave the lower bound $C \max\{n^{-4\alpha/(1+2\alpha)}, n^{-2\beta/(1+2\beta)}\}$ for the minimax risk. This bound is larger than the lower bound given in our Theorem 2 and is incorrect. This is due to a miscalculation on appendix C of their paper. A key step in that proof is to find some $d \geq 0$ such that

$$D = E\{[1 + \exp(\tfrac{1}{2}d + d^{1/2} N_1)]^{-1} (\tfrac{1}{2}d + d^{1/2} N_1)\} \neq 0.$$

In the above expression, $N_1$ denotes a standard normal random variable. But in fact

$$D = \int_{-\infty}^{\infty} \frac{(1/2)d + d^{1/2}x}{1 + \exp((1/2)d + d^{1/2}x)} \frac{1}{\sqrt{2\pi}} \exp\left(-\frac{x^2}{2}\right) dx$$

$$= \int_{-\infty}^{\infty} \frac{x}{1 + \exp(x)} \frac{1}{\sqrt{2\pi d}} \exp\left(-\frac{(x - (1/2)d)^2}{2d}\right) dx$$

$$= \int_{-\infty}^{\infty} \frac{x}{\exp(x/2) + \exp(-x/2)} \frac{1}{\sqrt{2\pi d}} \exp\left(-\frac{x^2}{2d} - \frac{d}{8}\right) dx.$$

This is an integral of an odd function which is identically 0 for all $d$.

**4. Discussion.** Variance function estimation in regression is more typically based on the residuals from a preliminary estimator $\hat{f}$ of the mean function. Such estimators have the form

(14)      $$\hat{V}(x) = \sum_i w_i(x)(y_i - \hat{f}(x_i))^2$$

where $w_i(x)$ are weight functions. A natural and common approach is to subtract in (14) an *optimal* estimator $\hat{f}$ of the mean function $f(x)$. See, for example, Hall and Carroll (1989), Neumann (1994), Ruppert, Wand, Holst and Hössjer (1997), and Fan and Yao (1998). When the unknown mean function is smooth, this approach often works well since the bias in $\hat{f}$ is negligible and $V$ can be estimated as well as when $f$ is identically zero. However, when the mean function is not smooth, using the residuals from an optimally smoothed $\hat{f}$ will lead to a sub-optimal estimator of $V$. For example, Hall and Carroll (1989) used a kernel estimator with optimal bandwidthfor $\hat{f}$ and showed that the resulting variance estimator attains the rate of

(15)      $$\max\{n^{-4\alpha/(2\alpha+1)}, n^{-2\beta/(2\beta+1)}\}$$

over $f \in \Lambda^\alpha(M_f)$ and $V \in \Lambda^\beta(M_V)$. This rate is strictly slower than the minimax rate when $\frac{4\alpha}{2\alpha+1} < \frac{2\beta}{2\beta+1}$ or equivalently, $\alpha < \frac{\beta}{2\beta+2}$.

Consider the example where $V$ belongs to a regular parametric family, such as $\{V(x) = \exp(ax + b) : a, b \in \mathbb{R}\}$. As Hall and Carroll have noted, this case is equivalent to the case of $\beta = \infty$ in results like Theorems 1 and 2. Then the rate of convergence for this estimator becomes nonparametric at $n^{-4\alpha/(2\alpha+1)}$ for $\alpha < 1/2$, while the optimal rate is the usual parametric rate $n^{-1/2}$ for all $\alpha \geq \frac{1}{4}$ and is $n^{-4\alpha}$ for $0 < \alpha < \frac{1}{4}$.

The main reason for the poor performance of such an estimator in the non-smooth setting is the "large" bias in $\hat{f}$. An optimal estimator $\hat{f}$ of $f$ balances the squared bias and variance. However, the bias and variance of $\hat{f}$ have significantly different effects on the estimation of $V$. The bias of $\hat{f}$ cannot be further reduced in the second stage smoothing of the squared residuals, while the variance of $\hat{f}$ can be incorporated easily. For $f \in \Lambda^\alpha(M_f)$, the maximum bias of an optimal estimator $\hat{f}$ is of order $n^{-\alpha/(2\alpha+1)}$ which becomes the dominant factor in the risk of $\hat{V}$ when $\alpha < \frac{\beta}{2\beta+2}$.

To minimize the effect of the mean function in such a setting one needs to use an estimator $\hat{f}(x_i)$ with minimal bias. Note that our approach is, in effect, using a very crude estimator $\hat{f}$ of $f$ with $\hat{f}(x_i) = y_{i+1}$. Such an estimator has high variance and low bias. As we have seen in Section 2, the large variance of $\hat{f}$ does not pose a problem (in terms of rates) for estimating $V$. Hence for estimating the variance function $V$ an optimal $\hat{f}$ is the one with minimum possible bias, not the one with minimum mean squared error. [Here we should of course exclude the obvious, and not useful, unbiased estimator $\hat{f}(x_i) = y_i$.]

Another implication of our results is that the unknown mean function does not have any first-order effect for estimating $V$ as long as $f$ has more than $1/4$ derivatives. When $\alpha > 1/4$, the variance estimator $\widehat{V}$ is essentially adaptive over $f \in \Lambda^\alpha(M_f)$ for all $\alpha > 1/4$. In other words, if $f$ is known to have more than $1/4$ derivatives, the variance function $V$ can be estimated with the same degree of first-order precision as if $f$ is completely known. However, when $\alpha < 1/4$, the rate of convergence for estimating $V$ is *entirely* determined by the degree of smoothness of the mean function $f$.

**5. Numerical results.** We now consider in this section the finite sample performance of our difference-based method for estimating the variance function. In particular we are interested in comparing the numerical performance of the difference-based estimator with the residual-based estimator of Fan and Yao (1998). The numerical results show that the performance of the difference-based estimator is somewhat inferior when the unknown mean function is very smooth. On the other hand, the difference-based estimator performs significantly better than the residual-based estimator when the mean function is not smooth.

Consider the model 1 where the variance function is $V(x) = (x - \frac{1}{2})^2 + \frac{1}{2}$ while there are four possible mean functions:

(i) $f_1(x) = 0$,
(ii) $f_2(x) = \frac{3}{4} * \sin(10\pi x)$,
(iii) $f_3(x) = \frac{3}{4} * \sin(20\pi x)$,
(iv) $f_4(x) = \frac{3}{4} * \sin(40\pi x)$.

The mean functions are arranged from a constant to much rougher sinusoid function; the "roughness" (the difficulty a particular mean function creates in estimation of the variance function $V$) is measured by the functional $R(f') = \int [f'(x)]^2 \, dx$ since the mean-related term in the asymptotic bias of the variance estimator $\hat{V}(x)$ is directly proportional to it. The numerical performance of the difference-based method had been investigated earlier in Levine (2006) for a slightly different set of mean functions.

TABLE 1
*Performance under the changing curvature of the mean function*

| | | Median CDMSE | |
|---|---|---|---|
| **Mean function** | $R(f')$ | **Fan–Yao method** | **Our method** |
| $f = 0$ | 0 | 0.00299 | 0.00376 |
| $f = \frac{3}{4}\sin(10\pi x)$ | 278.15 | 0.07161 | 0.00344 |
| $f = \frac{3}{4}\sin(20\pi x)$ | 1110.89 | 0.08435 | 0.00384 |
| $f = \frac{3}{4}\sin(40\pi x)$ | 4441.88 | 0.08363 | 0.00348 |

For comparison purposes, the same four combinations of the mean and variance functions are investigated using the two-step method described in Fan and Yao (1998). We expect this method to perform better than the difference-based method in the case of a constant mean function, but to get progressively worse as the roughness of the mean function considered increases. The following table summarizes results of simulations using both methods. In this case, the bandwidths for estimating the mean and variance functions were selected using a $K$-fold cross-validation with $K = 10$. We consider the fixed equidistant design $x_i = \frac{i}{n}$ on $[0, 1]$ where the sample size is $n = 1000$; 100 simulations are performed and the bandwidth $h$ is selected using a $K$-fold cross-validation with $K = 10$. The performance of both methods is measured using the cross-validation discrete mean squared error (CDMSE) that is defined as

$$\text{(16)} \quad \text{CDMSE} = n^{-1} \sum_{i=1}^{n} [\hat{V}_{h_{\text{CV}}}(x_i) - V(x_i)]^2$$

with $h_{\text{CV}}$ being the $K$-fold cross-validation bandwidth. We report the median CDMSE for variance function estimators based on 100 simulations. Table 1 provides the summary of the performance.

It is easily seen from the table that the two-step method of Fan and Yao, based on estimating the variance using squared residuals from the mean function estimation, tends to perform slightly better when the mean function is very smooth but noticeably worse when it is rougher. Note that here we only use the first-order differences. The performance of the difference based estimator can be improved in the case of smooth mean function by using higher order differences. The Fan–Yao method performs about 26% better in the first case of the constant mean function. However, the risk (CDMSE) of the difference based method is over 95% smaller than the risk of the Fan–Yao method for the second mean function. In the rougher cases, the difference is approximately the same. The CDMSE of the difference based method is over 95% and 96% less than the corresponding risk of the residual based method for the third and fourth mean functions, respectively.

## 6. Proofs.

6.1. *Upper bound*: *Proof of Theorem* 1. We shall only prove (8). Inequality (9) is a direct consequence of (8). Recall that

$$D_i^2 = \delta_i^2 + 2V_i + 2V_i(\epsilon_i^2 - 1) + 2\sqrt{2}\delta_i V_i^{1/2}\epsilon_i,$$

where $\delta_i = f(x_i) - f(x_{i+1})$, $V_i^{1/2} = \sqrt{1/2(V(x_i) + V(x_{i+1}))}$ and

$$\epsilon_i = (V(x_i) + V(x_{i+1}))^{-1/2}(V^{1/2}(x_i)z_i - V^{1/2}(x_{i+1})z_{i+1}).$$

Without loss of generality, suppose $h = n^{-1/(1+2\beta)}$. It is easy to see that for any $x_* \in [0,1]$, $\sum_i K_i^h(x_*) = 1$, and when $x_* \geq (x_i + x_{i+1})/2 + h$ or $x_* \leq (x_i + x_{i-1})/2 - h$, $K_i^h(x_*)$ equals 0. Suppose $k < \widehat{k}$, we also have

$$\left(\sum_i |K_i^h(x_*)|\right)^2 \leq 2nh \sum_i (K_i^h(x_*))^2$$

$$\leq 2\int_{-1}^1 K_*^2(u)\,du$$

$$\leq 2\widehat{k},$$

where $K_*(u) = K(u)$ when $x_* \in (h, 1-h)$; $K_*(u) = K_t(u)$ when $x_* = th$ for some $t \in [0,1]$; and $K_*(u) = K_t(-u)$ when $x_* = 1 - th$ for some $t \in [0,1]$.

The second inequality above is obtained as follows. For the sake of simplicity, assume that $K_* = K$; the same argument can be repeated for boundary kernels as well. Using the definition of $K_i^h(x_*)$, we note that it can be rewritten as $\int_{(x_i+x_{i-1})/2}^{(x_i+x_{i+1})/2} \frac{1}{nh} K(\frac{x-u}{h})\,d(nu)$. Since the last integral is taken with respect to the probability measure on the interval $[\frac{x_i+x_{i-1}}{2}, \frac{x_i+x_{i+1}}{2}]$, we can apply Jensen's inequality to obtain

$$(K_i^h(x_*))^2 \leq \frac{1}{nh}\int_{(x_i+x_{i-1})/2}^{(x_i+x_{i+1})/2} K^2\left(\frac{x-u}{h}\right)d(nu)$$

$$= \frac{1}{(nh)^2}\int_{(x_i+x_{i-1})/2}^{(x_i+x_{i+1})/2} K^2\left(\frac{x-u}{h}\right)du.$$

Thus,

$$\left(\sum_i |K_i^h(x_*)|\right)^2 \leq \frac{2}{h}\sum_i \int_{(x_i+x_{i-1})/2}^{(x_i+x_{i+1})/2} K^2\left(\frac{x-u}{h}\right)du$$

$$= 2 \int_{-1}^{1} K^2(u)\, du.$$

For all $f \in \Lambda^\alpha(M_f)$ and $V \in \Lambda^\beta(M_V)$, the mean squared error of $\hat{V}$ at $x_*$ satisfies

$$E(\widehat{V}(x_*) - V(x_*))^2$$

$$= E\left(\sum_{i=1}^{n-1} K_i^h(x_*)(\tfrac{1}{2}D_i^2 - V(x_*))\right)^2$$

$$= E\Bigg\{\sum_{i=1}^{n-1} K_i^h(x_*)\tfrac{1}{2}\delta_i^2 + \sum_{i=1}^{n-1} K_i^h(x_*)(V_i - V(x_*))$$

$$+ \sum_{i=1}^{n-1} K_i^h(x_*)V_i(\epsilon_i^2 - 1) + \sum_{i=1}^{n-1} K_i^h(x_*)\sqrt{2}\delta_i V_i^{1/2}\epsilon_i\Bigg\}^2$$

$$\leq 4\left(\sum_{i=1}^{n-1} K_i^h(x_*)\tfrac{1}{2}\delta_i^2\right)^2 + 4\left(\sum_{i=1}^{n-1} K_i^h(x_*)(V_i - V(x_*))\right)^2$$

$$+ 4E\left(\sum_{i=1}^{n-1} K_i^h(x_*)V_i(\epsilon_i^2 - 1)\right)^2 + 4E\left(\sum_{i=1}^{n-1} K_i^h(x_*)\sqrt{2}\delta_i V_i^{1/2}\epsilon_i\right)^2.$$

Suppose $\alpha \leq 1/4$, otherwise $n^{-4\alpha} < n^{-2\beta/(1+2\beta)}$ for any $\beta$. Since for any $i$, $|\delta_i| = |f(x_i) - f(x_{i+1})| \leq M_f |x_i - x_{i+1}|^\alpha = M_f n^{-\alpha}$, we have

$$4\left(\sum_{i=1}^{n-1} K_i^h(x_*)\tfrac{1}{2}\delta_i^2\right)^2 \leq 4\left(\sum_{i=1}^{n-1} |K_i^h(x_*)|\tfrac{1}{2}M_f^2 n^{-2\alpha}\right)^2 \leq 2\widehat{k} M_f^4 n^{-4\alpha}.$$

Note that for any $x, y \in [0, 1]$, Taylor's theorem yields

$$\left| V(x) - V(y) - \sum_{j=1}^{\lfloor \beta \rfloor} \frac{V^{(j)}(y)}{j!}(x-y)^j \right|$$

$$= \left| \int_y^x \frac{(x-u)^{\lfloor \beta \rfloor - 1}}{(\lfloor \beta \rfloor - 1)!} (V^{(\lfloor \beta \rfloor)}(u) - V^{(\lfloor \beta \rfloor)}(y))\, du \right|$$

$$\leq \left| \int_y^x \frac{(x-u)^{\lfloor \beta \rfloor - 1}}{(\lfloor \beta \rfloor - 1)!} M_V |x-y|^{\beta - \lfloor \beta \rfloor}\, du \right|$$

$$\leq \frac{M_V}{\lfloor \beta \rfloor!} |x-y|^\beta.$$

So,

$$V_i - V(x_*) = \frac{1}{2}\left(V\left(\frac{i}{n}\right) + V\left(\frac{i+1}{n}\right)\right) - V(x_*)$$

$$\leq \frac{1}{2}\sum_{j=1}^{\lfloor \beta \rfloor} \frac{V^{(j)}(x_*)}{j}\left(\left(\frac{i}{n}-x_*\right)^j + \left(\frac{i+1}{n}-x_*\right)^j\right)$$

$$+\frac{1}{2}M_V\left|\frac{i}{n}-x_*\right|^\beta + \frac{1}{2}M_V\left|\frac{i+1}{n}-x_*\right|^\beta$$

and

$$V_i - V(x_*) \geq \frac{1}{2}\sum_{j=1}^{\lfloor \beta \rfloor} \frac{V^{(j)}(x_*)}{j}\left(\left(\frac{i}{n}-x_*\right)^j + \left(\frac{i+1}{n}-x_*\right)^j\right)$$

$$-\frac{1}{2}M_V\left|\frac{i}{n}-x_*\right|^\beta - \frac{1}{2}M_V\left|\frac{i+1}{n}-x_*\right|^\beta.$$

Since the kernel functions have vanishing moments, for $j = 1, 2, \ldots, \lfloor \beta \rfloor$, when $n$ large enough

$$\left|\sum_{i=1}^{n-1} K_i^h(x_*)\left(\frac{i}{n}-x_*\right)^j\right|$$

$$= \left|\sum_{i=1}^{n-1}\int_{(x_i+x_{i-1})/2}^{(x_i+x_{i+1})/2} \frac{1}{h}K\left(\frac{x_*-u}{h}\right)\left(\frac{i}{n}-x_*\right)^j du\right|$$

$$= \left|\int_0^1 \frac{1}{h}K\left(\frac{x_*-u}{h}\right)(u-x_*)^j du\right.$$

$$\left. + \sum_{i=1}^{n-1}\int_{(x_i+x_{i-1})/2}^{(x_i+x_{i+1})/2} \frac{1}{h}K\left(\frac{x_*-u}{h}\right)\left[\left(\frac{i}{n}-x_*\right)^j - (u-x_*)^j\right] du\right|$$

$$= \left|\sum_{i=1}^{n-1}\int_{(x_i+x_{i-1})/2}^{(x_i+x_{i+1})/2} \frac{1}{h}K\left(\frac{x_*-u}{h}\right)\left[\left(\frac{i}{n}-x_*\right)^j - (u-x_*)^j\right] du\right|$$

$$\leq c'\sum_{i=1}^{n-1}\int_{(x_i+x_{i-1})/2}^{(x_i+x_{i+1})/2} \left|\frac{1}{h}\left(\frac{x_*-u}{h}\right)\right| \times \frac{j}{n} du = c'n^{-1}$$

for some generic constant $c' > 0$. Similarly, $\sum_{i=1}^{n-1} K_i^h(x_*)(\frac{i+1}{n}-x_*)^j \leq c'n^{-1}$. So,

$$\left|\sum_{i=1}^{n-1} K_i^h(x_*)\left(\sum_{j=1}^{\lfloor \beta \rfloor} \frac{V^{(j)}(x_*)}{j}\left(\left(\frac{i}{n}-x_*\right)^j + \left(\frac{i+1}{n}-x_*\right)^j\right)\right)\right| \leq \widehat{C}n^{-1}$$

for some constant $\widehat{C} > 0$ which does not depend on $x_*$. Note that $V^{\lfloor \beta \rfloor}$ satisfies Hölder condition with exponent $0 < \alpha' = \alpha - \lfloor \alpha \rfloor < 1$ and is, therefore, continuous on $[0,1]$ and bounded. Then we have

$$4\left(\sum_{i=1}^{n-1} K_i^h(x_*)(V_i - V(x_*))\right)^2$$

$$\leq 2\widehat{C}^2 n^{-2} + 2M_V^2 \left(\sum_{i=\lfloor n(x_*-h) \rfloor}^{\lfloor n(x_*+h) \rfloor + 1} |K_i^h(x_*)| \left(\left|\frac{i}{n} - x_*\right|^\beta + \left|\frac{i+1}{n} - x_*\right|^\beta\right)\right)^2$$

$$\leq 2\widehat{C}^2 n^{-2} + 2M_V^2 \left(\sum_{i=\lfloor n(x_*-h) \rfloor}^{\lfloor n(x_*+h) \rfloor + 1} |K_i^h(x_*)| \left(\left|h + \frac{1}{n}\right|^\beta + \left|h + \frac{2}{n}\right|^\beta\right)\right)^2$$

$$\leq 2\widehat{C}^2 n^{-2} + 8 \times 3^{2\beta} M_V^2 n^{-2\beta/(1+2\beta)} \times (2\widehat{k}).$$

The last inequality is due to the fact $0 < h + \frac{1}{n} < h + \frac{2}{n} < 3h$. On the other hand, notice that $\epsilon_1, \epsilon_3, \epsilon_5, \ldots$ are independent and $\epsilon_2, \epsilon_4, \epsilon_6, \ldots$ are independent, we have

$$4E\left(\sum_{i=1}^{n-1} K_i^h(x_*)\sqrt{2}\delta_i V_i^{1/2}\epsilon_i\right)^2 = 4\operatorname{Var}\left(\sum_{i=1}^{n-1} K_i^h(x_*)\sqrt{2}\delta_i V_i^{1/2}\epsilon_i\right)$$

$$\leq 16 \sum_{i=\lfloor n(x_*-h) \rfloor}^{\lfloor n(x_*+h) \rfloor + 1} (K_i^h(x_*))^2 \delta_i^2 V_i$$

$$\leq 16 M_f^2 M_V n^{-2\alpha - 2\beta/(1+2\beta)} \times \widehat{k}$$

and

$$4E\left(\sum_{i=1}^{n-1} K_i^h(x_*)V_i(\epsilon_i^2 - 1)\right)^2 = 4\operatorname{Var}\left(\sum_{i=1}^{n-1} K_i^h(x_*)V_i(\epsilon_i^2 - 1)\right)$$

$$\leq 8M_V^2 \mu_4 \sum_{i=1}^{n-1} (K_i^h(x_*))^2$$

$$\leq 8M_V^2 \mu_4 \frac{1}{nh}\widehat{k}$$

$$= 8M_V^2 \mu_4 n^{-2\beta/(1+2\beta)} \times \widehat{k}$$

where $\mu_4$ denotes the uniform bound for the fourth moments of the $\epsilon_i$.

Putting the four terms together we have, uniformly for all $x_* \in [0,1]$, $f \in \Lambda^\alpha(M_f)$ and $V \in \Lambda^\beta(M_V)$,

$$E(\widehat{V}(x_*) - V(x_*))^2$$

$$\leq 2\widehat{k}M_f^4 n^{-4\alpha} + 2\widehat{C}^2 n^{-2} + 8 \times 3^{2\beta} M_V^2 n^{-2\beta/(1+2\beta)} \times (2\widehat{k})$$

$$+ 8M_V^2\mu_4 n^{-2\beta/(1+2\beta)}\widehat{k} + 16M_f^2 M_V n^{-2\alpha-2\beta/(1+2\beta)}\widehat{k}$$
$$= C_0 \cdot \max\{n^{-4\alpha}, n^{-2\beta/(1+2\beta)}\}$$

for some constant $C_0 > 0$. This proves (8).

6.2. *Lower bound*: *Proof of Theorem* 2. We shall only prove the lower bound for the pointwise squared error loss. The same proof with minor modifications immediately yields the lower bound under the integrated squared error. Note that, to prove inequality (13), we only need to focus on the case where $\alpha < 1/4$, otherwise $n^{-2\beta/(1+2\beta)}$ is always greater than $n^{-4\alpha}$ for sufficiently large $n$ and then (13) follows directly from (12).

For a given $0 < \alpha < 1/4$, there exists an integer $q$ such that $(q+1)\alpha > 1$. For convenience we take $q$ to be an odd integer. From Lemma 1, there is a positive constant $B < \infty$ and a symmetric distribution $G$ on $[-B, B]$ such that $G$ and $N(0,1)$ have the same first $q$ moments. Let $r_i$, $i = 1,\ldots,n$, be independent variables with the distribution $G$. Set $\theta_n = \frac{M_f}{2B}n^{-\alpha}$, $f_0 \equiv 0$, $V_0(x) \equiv 1 + \theta_n^2$ and $V_1(x) \equiv 1$. Let $g(x) = 1 - 2n|x|$ for $x \in [-\frac{1}{2n}, \frac{1}{2n}]$ and 0 otherwise. Define the random function $f_1$ by

$$f_1(x) = \sum_{i=1}^{n} \theta_n r_i g(x - x_i) I(0 \leq x \leq 1).$$

Then it is easy to see that $f_1$ is in $\Lambda^\alpha(M_f)$ for all realizations of $r_i$. Moreover, $f_1(x_i) = \theta_n r_i$ are independent and identically distributed.

Now consider testing the following hypotheses:

$$H_0 : y_i = f_0(x_i) + V_0^{1/2}(x_i)\nu_i, \qquad i = 1,\ldots,n,$$
$$H_1 : y_i = f_1(x_i) + V_1^{1/2}(x_i)\nu_i, \qquad i = 1,\ldots,n,$$

where $\nu_i$ are independent $N(0,1)$ variables which are also independent of the $r_i$'s. Denote by $P_0$ and $P_1$ the joint distributions of $y_i$'s under $H_0$ and $H_1$, respectively. Note that for any estimator $\widehat{V}$ of $V$,

$$\max\{E(\widehat{V}(x_*) - V_0(x_*))^2, E(\widehat{V}(x_*) - V_1(x_*))^2\}$$
(17)
$$\geq \frac{1}{16}\rho^4(P_0, P_1)(V_0(x_*) - V_1(x_*))^2$$
$$= \frac{1}{16}\rho^4(P_0, P_1)\frac{M_f^4}{16B^4}n^{-4\alpha}$$

where $\rho(P_0, P_1)$ is the Hellinger affinity between $P_0$ and $P_1$. See, for example, Le Cam (1986). Let $p_0$ and $p_1$ be the probability density function of $P_0$ and $P_1$ with respect to the Lebesgue measure $\mu$, then $\rho(P_0, P_1) = \int \sqrt{p_0 p_1}\, d\mu$. The minimax lower bound (13) follows immediately from the two-point bound (17) if we show that for any $n$, the Hellinger affinity $\rho(P_0, P_1) \geq C$ for some constant $C > 0$. ($C$ may depend on $q$, but does not depend on $n$.)

Note that under $H_0$, $y_i \sim N(0, 1 + \theta_n^2)$ and its density $d_0$ can be written as

$$d_0(t) \triangleq \frac{1}{\sqrt{1+\theta_n^2}} \varphi\left(\frac{t}{\sqrt{1+\theta_n^2}}\right) = \int \varphi(t - v\theta_n)\varphi(v)\, dv.$$

Under $H_1$, the density of $y_i$ is $d_1(t) \triangleq \int \varphi(t - v\theta_n) G(dv)$.

It is easy to see that $\rho(P_0, P_1) = (\int \sqrt{d_0 d_1}\, d\mu)^n$, since the $y_i$'s are independent variables. Note that the Hellinger affinity is bounded below by the total variation affinity

$$\int \sqrt{d_0(t) d_1(t)}\, dt \geq 1 - \tfrac{1}{2} \int |d_0(t) - d_1(t)|\, dt.$$

Taylor's expansion yields

$$\varphi(t - v\theta_n) = \varphi(t) \left( \sum_{k=0}^{\infty} v^k \theta_n^k \frac{H_k(t)}{k!} \right),$$

where $H_k(t)$ is the corresponding Hermite polynomial. And from the construction of the distribution $G$,

$$\int v^i G(dv) = \int v^i \varphi(v)\, dv \qquad \text{for } i = 0, 1, \ldots, q.$$

So,

$$|d_0(t) - d_1(t)|$$
$$= \left| \int \varphi(t - v\theta_n) G(dv) - \int \varphi(t - v\theta_n)\varphi(v)\, dv \right|$$
(18)
$$= \left| \int \varphi(t) \sum_{i=0}^{\infty} \frac{H_i(t)}{i!} v^i \theta_n^i G(dv) - \int \varphi(t) \sum_{i=0}^{\infty} \frac{H_i(t)}{i!} v^i \theta_n^i \varphi(v)\, dv \right|$$
$$= \left| \int \varphi(t) \sum_{i=q+1}^{\infty} \frac{H_i(t)}{i!} v^i \theta_n^i G(dv) - \int \varphi(t) \sum_{i=q+1}^{\infty} \frac{H_i(t)}{i!} v^i \theta_n^i \varphi(v)\, dv \right|$$
$$\leq \left| \int \varphi(t) \sum_{i=q+1}^{\infty} \frac{H_i(t)}{i!} v^i \theta_n^i G(dv) \right| + \left| \int \varphi(t) \sum_{i=q+1}^{\infty} \frac{H_i(t)}{i!} v^i \theta_n^i \varphi(v)\, dv \right|.$$

Suppose $q + 1 = 2p$ for some integer $p$, it can be seen that

$$\left| \int \varphi(t) \sum_{i=q+1}^{\infty} \frac{H_i(t)}{i!} v^i \theta_n^i G(dv) \right| = \left| \int \varphi(t) \sum_{i=p}^{\infty} \frac{H_{2i}(t)}{(2i)!} \theta_n^{2i} v^{2i} G(dv) \right|$$
$$\leq \varphi(t) \sum_{i=p}^{\infty} \left| \frac{H_{2i}(t)}{(2i)!} \theta_n^{2i} \right| \left| \int v^{2i} G(dv) \right|$$

$$\leq \varphi(t) \sum_{i=p}^{\infty} \left| \frac{H_{2i}(t)}{(2i)!} \right| \theta_n^{2i} B^{2i}$$

and

$$\left| \int \varphi(t) \sum_{i=q+1}^{\infty} \frac{H_i(t)}{i!} v^i \theta_n^i \varphi(v) \, dv \right| = \left| \int \varphi(t) \sum_{i=p}^{\infty} \frac{H_{2i}(t)}{(2i)!} \theta_n^{2i} v^{2i} \varphi(v) \, dv \right|$$

$$\leq \varphi(t) \sum_{i=p}^{\infty} \left| \frac{H_{2i}(t)}{(2i)!} \theta_n^{2i} \right| \left| \int v^{2i} \varphi(v) \, dv \right|$$

$$= \left| \varphi(t) \sum_{i=p}^{\infty} H_{2i}(t) \theta_n^{2i} \frac{1}{2^i \cdot i!} \right|$$

$$\leq \varphi(t) \sum_{i=p}^{\infty} \left| \frac{H_{2i}(t)}{2^i \cdot i!} \right| \theta_n^{2i}.$$

So from (18),

$$|d_0(t) - d_1(t)| \leq \varphi(t) \sum_{i=p}^{\infty} \left| \frac{H_{2i}(t)}{(2i)!} \right| \theta_n^{2i} B^{2i} + \varphi(t) \sum_{i=p}^{\infty} \left| \frac{H_{2i}(t)}{2^i \cdot i!} \right| \theta_n^{2i}$$

and then

$$\int \sqrt{d_0(t) d_1(t)} \, dt$$

(19)
$$\geq 1 - \frac{1}{2} \int \left( \varphi(t) \sum_{i=p}^{\infty} \left| \frac{H_{2i}(t)}{(2i)!} \right| \theta_n^{2i} B^{2i} + \varphi(t) \sum_{i=p}^{\infty} \left| \frac{H_{2i}(t)}{2^i \cdot i!} \right| \theta_n^{2i} \right) dt$$

$$= 1 - \frac{1}{2} \int \varphi(t) \sum_{i=p}^{\infty} \left| \frac{H_{2i}(t)}{(2i)!} \right| \theta_n^{2i} B^{2i} \, dt - \frac{1}{2} \int \varphi(t) \sum_{i=p}^{\infty} \left| \frac{H_{2i}(t)}{2^i \cdot i!} \right| \theta_n^{2i} \, dt.$$

Since $\int t^{2i} \phi(t) \, dt = (2i-1)!!$ where $(2i-1)!! \triangleq (2i-1) \times (2i-3) \times \cdots \times 3 \times 1$, for the Hermite polynomial $H_{2i}$ we have

$$\int \varphi(t) |H_{2i}(t)| \, dt$$

$$= \int \varphi(t) \left| (2i-1)!! \times \left[ 1 + \sum_{k=1}^{i} \frac{(-2)^k i(i-1)\cdots(i-k+1)}{(2k)!} t^{2k} \right] \right| dt$$

$$\leq \int \varphi(t) \left[ (2i-1)!! \times \left( 1 + \sum_{k=1}^{i} \frac{2^k i(i-1)\cdots(i-k+1)}{(2k)!} t^{2k} \right) \right] dt$$

$$= (2i-1)!! \times \left( 1 + \sum_{k=1}^{i} \frac{2^k i(i-1)\cdots(i-k+1)}{(2k)!} \int t^{2k} \varphi(t) \, dt \right)$$

$$= (2i-1)!! \times \left(1 + \sum_{k=1}^{i} \frac{2^k i(i-1)\cdots(i-k+1)}{(2k)!}(2k-1)!!\right)$$

$$= (2i-1)!! \times \left(1 + \sum_{k=1}^{i} \frac{i(i-1)\cdots(i-k+1)}{k!}\right)$$

$$= 2^i \times (2i-1)!!.$$

For sufficiently large $n$, $\theta_n < 1/2$ and it then follows from the above inequality that

$$\int \varphi(t) \sum_{i=p}^{\infty} \left|\frac{H_{2i}(t)}{(2i)!}\right| \theta_n^{2i} B^{2i} \, dt \leq \sum_{i=p}^{\infty} \frac{\theta_n^{2i} B^{2i}}{(2i)!} \int \varphi(t) |H_{2i}(t)| \, dt$$

$$\leq \sum_{i=p}^{\infty} \frac{\theta_n^{2i} B^{2i}}{(2i)!} 2^i \times (2i-1)!!$$

$$= \theta_n^{2p} \sum_{i=p}^{\infty} \frac{B^{2i} \theta_n^{2i-2p}}{i!} \leq \theta_n^{2p} \times e^{B^2}$$

and

$$\int \varphi(t) \sum_{i=p}^{\infty} \left|\frac{H_{2i}(t)}{2^i \cdot i!}\right| \theta_n^{2i} \, dt \leq \sum_{i=p}^{\infty} \frac{\theta_n^{2i}}{2^i \cdot i!} \int \varphi(t) |H_{2i}(t)| \, dt$$

$$\leq \sum_{i=p}^{\infty} \frac{\theta_n^{2i}}{2^i \cdot i!} 2^i \times (2i-1)!!$$

$$= \theta_n^{2p} \sum_{i=p}^{\infty} \frac{(2i-1)!!}{i!} \theta_n^{2i-2p}$$

$$\leq \theta_n^{2p} \sum_{i=p}^{\infty} 2^i \times \theta_n^{2i-2p} \leq \theta_n^{2p} \sum_{i=p}^{\infty} 2^i \times \left(\frac{1}{2}\right)^{2i-2p}$$

$$= \theta_n^{2p} \times 2^{2p+1}.$$

Then from (19)

$$\int \sqrt{d_0(t)d_1(t)} \, dt \geq 1 - \theta_n^{2p}(\tfrac{1}{2}e^{B^2} + 2^{2p}) \triangleq 1 - c\theta_n^{q+1},$$

where $c$ is a constant that only depends on $q$. So

$$\rho(P_0, P_1) = \left(\int \sqrt{d_0(t)d_1(t)} \, dt\right)^n \geq (1 - c\theta_n^{q+1})^n = (1 - cn^{-\alpha(q+1)})^n.$$

Since $\alpha(q+1) \geq 1$, $\lim_{\to \infty}(1 - cn^{-\alpha(q+1)})^n \geq e^{-c} > 0$ and the theorem then follows from (17).

**Acknowledgments.** We thank the Editor and three referees for thorough and useful comments which have helped to improve the presentation of the paper. We would like to thank William Studden for discussions on the finite moment matching problem and for the reference Karlin and Studden (1966).

## REFERENCES

Brown, L. D. and Levine, M. (2006). Variance estimation in nonparametric regression via the difference sequence method. *Ann. Statist.* To appear.

Fan, J. and Gijbels, I. (1996). *Local Polynomial Modelling and Its Applications.* Chapman and Hall, London. MR1383587

Fan, J. and Yao, Q. (1998). Efficient estimation of conditional variance functions in stochastic regression. *Biometrika* **85** 645–660. MR1665822

Gasser, T. and Müller, H. G. (1979). Kernel estimation of regression functions. *Smoothing Techniques for Curve Estimation. Lecture Notes in Math.* **757** 23–68. Springer, New York. MR0564251

Gasser, T. and Müller, H. G. (1984). Estimating regression functions and their derivatives by the kernel method. *Scand. J. Statist.* **11** 197–211. MR0767241

Gasser, T., Müller, H. G. and Mammitzsch, V. (1985). Kernels for nonparametric curve estimation. *J. Roy. Statist. Soc. B* **47** 238–252. MR0816088

Hall, P. and Carroll, R. J. (1989). Variance function estimation in regression: The effect of estimating the mean. *J. Roy. Statist. Soc. Ser. B* **51** 3–14. MR0984989

Hall, P., Kay, J. and Titterington, D. M. (1990). Asymptotically optimal difference-based estimation of variance in nonparametric regression. *Biometrika* **77** 521–528. MR1087842

Karlin, S. and Studden, W. J. (1966). *Tchebycheff Systems*: *With Applications in Analysis and Statistics.* Interscience, New York. MR0204922

Le Cam, L. (1986). *Asymptotic Methods in Statistical Decision Theory.* Springer, New York. MR0856411

Levine, M. (2006). Bandwidth selection for a class of difference-based variance estimators in the nonparametric regression: A possible approach. *Comput. Statist. Data Anal.* **50** 3405–3431. MR2236857

Müller, H. G. and Stadtmüller, U. (1987). Estimation of heteroscedasticity in regression analysis. *Ann. Statist.* **15** 610–625. MR0888429

Müller, H. G. and Stadtmüller, U. (1993). On variance function estimation with quadratic forms. *J. Statist. Plann. Inference* **35** 213–231. MR1220417

Munk, A., Bissantz, N., Wagner, T. and Freitag, G. (2005). On difference based variance estimation in nonparametric regression when the covariate is high dimensional. *J. Roy. Statist. Soc. Ser. B* **67** 19–41. MR2136637

Neumann, M.(1994). Fully data-driven nonparametric variance estimators. *Statistics* **25** 189–212. MR1366825

Rice, J. (1984). Bandwidth choice for nonparametric kernel regression. *Ann. Statist.* **12** 1215–1230. MR0760684

Ruppert, D., Wand, M. P., Holst, U. and Hössjer, O. (1997). Local polynomial variance function estimation. *Technometrics* **39** 262–273. MR1462587

von Neumann, J.(1941). Distribution of the ratio of the mean squared successive difference to the variance. *Ann. Math. Statist.* **12** 367–395. MR0006656

von Neumann, J. (1942). A further remark concerning the distribution of the ratio of the mean squared successive difference to the variance. *Ann. Math. Statist.* **13** 86–88. MR0006657

L. Wang  
L. D. Brown  
T. T. Cai  
Department of Statistics  
The Wharton School  
University of Pennsylvania  
Philadelphia, Pennsylvania 19104  
USA  
E-mail: tcai@wharton.upenn.edu

M. Levine  
Department of Statistics  
Purdue University  
West Lafayette, Indiana 47907  
USA  
E-mail: mlevins@stat.purdue.edu